\documentclass[12pt]{article}%
\usepackage{amssymb}
\usepackage{amsfonts}
\usepackage{amsmath}
\usepackage{graphicx}%
\providecommand{\U}[1]{\protect\rule{.1in}{.1in}}
\newtheorem{theorem}{Theorem}[section]

\newtheorem{conjecture}[theorem]{Conjecture}

\newtheorem{problem}[theorem]{Problem}
\newtheorem{question}{Question}
\newtheorem{proposition}[theorem]{Proposition}

\begin{document}

\title{\textbf{The Eternal Game Chromatic Number of a Graph}}
\author{William F. Klostermeyer\\
      School of Computing\\
      University of North Florida \\ Jacksonville, FL 32224-2669\\{\small klostermeyer@hotmail.com}
    \and Hannah Mendoza\\Wake Forest University\\Winston Salem, NC 27109\\{\small mendhm17@wfu.edu}}

\maketitle

\begin{abstract} Game coloring is a well-studied two-player game in which each player properly colors one vertex of a graph at a time until all the vertices are colored. An ``eternal'' version of game coloring is introduced in this paper in which the vertices are colored and re-colored from a color set over a sequence of rounds. In a given round, each vertex is colored, or re-colored, once, so that a proper coloring is maintained. Player 1 wants to maintain a  proper coloring forever, while player 2 wants to  force the coloring process to fail. The eternal game chromatic number of a graph $G$ is defined to be the minimum number of colors needed in the color set so that player 1 can always win the game on $G$. We consider several variations of this new game and show its behavior on some elementary classes of graphs.
\end{abstract}



\section{Introduction}

The {\it{eternal graph coloring problem}} was introduced by Klostermeyer in \cite{WFK} and is defined as follows. Let $G=(V, E)$ be a finite, undirected graph with $n$ vertices. Let $f_0: V \rightarrow Z^+$ be a proper vertex coloring of $G$. An infinite sequence of vertex {\it requests} $R=r_1, r_2, \ldots$ must be handled as follows. After $r_i$ is revealed, the color assigned to $r_i$, $f_{i-1}(r_i)$, must be changed, yielding a proper coloring $f_i$. We may formalize this problem as a two-player game: player 1 chooses an initial proper coloring of $G$ with at most $k$ colors; the players then alternate turns with player 2  (the {\it adversary}) choosing a vertex and player 1 changing the color of that vertex while still maintaining a proper coloring of $G$ with at most $k$ colors. Player 2 wins if player 1 has no valid move at some turn and player 1 wins otherwise.  The smallest positive integer $k$ for which player 1 can win the game on graph $G$, for any possible sequence of moves by player 2, is the {\it eternal chromatic number} of $G$, is denoted $\chi^\infty(G)$.  In other words, how many colors are needed to ensure player 1 can win the game?
The game may be viewed as starting from some initial proper coloring, with player 1 maintaining $\chi^\infty(G)$ independent sets, dynamically moving vertices from one set to another as requests from player 2 are made so that each set remains an independent set at all times.

The {\it vertex coloring game} was introduced in 1981 by Brams (see \cite{Gard}) and rediscovered later by Bodlaender \cite{Bod}. It is a two-player game played according to the following rules:\\

\noindent $\bullet$ Alice and Bob properly color the vertices of a graph $G$ with a fixed set of $k$ colors (the colors are assumed to be integers in the range $1, \ldots, k$).\\
\noindent $\bullet$  Alice and Bob take turns, coloring properly any uncolored vertex (in the standard version of the game, Alice begins).\\
\noindent $\bullet$ If a vertex $v$ is impossible to color properly (that is, for any color in the color set, $v$ has a neighbor colored with that color), then Bob wins.\\
\noindent $\bullet$ If all vertices in the graph are colored properly, then Alice wins.\\

The game chromatic number of $G$, denoted by $\chi _{g}(G)$, is the minimum number of colors needed in the color set to guarantee that Alice wins the vertex coloring game on $G$. 

The purpose of this paper is to define a new graph coloring game which extends game coloring to an ``eternal'' model of graph coloring that is played by two players over a series of rounds. We define several versions of this game, present some basic results, and pose some questions.

\section{Eternal Game Chromatic Number}

We define the {\it eternal vertex coloring game} on finite, undirected graph $G$ as follows.

\noindent $\bullet$ Alice and Bob properly color the vertices of a graph G with a fixed set of $k$ colors (the colors are assumed to be integers in the range $1, \ldots, k$).\\
\noindent $\bullet$ Alice and Bob take turns, choosing and coloring (or re-coloring) properly one vertex at a time, with Alice going first:\\
\indent $\bullet$ if a vertex $v$ is uncolored and is chosen by Alice or Bob, then it must be properly colored\\
\indent $\bullet$ if a vertex $v$ is colored color $c$ and is chosen by Alice or Bob, then it must be assigned a color other than $c$ which is not used on one of its neighbors (this new color is chosen by the player whose turn it is) \\
\indent $\bullet$ all vertices must be chosen  $t$ times prior to any vertex being chosen $t+1$ times, for every integer $t \geq 1$\\
\noindent $\bullet$ If a vertex $v$ is chosen and is impossible to color properly (that is, for every color in the color set other than its current color, $v$ has a neighbor colored with that color), then Bob wins.\\
\noindent  $\bullet$ If Bob cannot win the game, then Alice wins.\\

In other words, the colorings take place over a series of {\it rounds}: each vertex is colored, or re-colored, once during each round. The number of rounds is infinite if Alice wins. We note that if $n$ is odd, the second round, and all even numbered rounds, have Bob going first. We shall sometimes refer to this as the A-game.

Let $\chi_g^\infty(G)$ be the minimum number of colors needed in the initial color set to guarantee Alice to win the eternal vertex coloring game on $G$. We call  $\chi_g^\infty$ the {\it eternal game chromatic number}.  As an example,  trivially $\chi_g^\infty(K_2)=3$. Not so trivially, $\chi_g^\infty(P_3)=3$. To see this,
let $P_3=uvw$. Alice colors $v$ red, Bob colors $u$ blue, and Alice colors $w$ blue. There are now two cases for Bob's second turn:

Case 1) If Bob colors $v$ green, then Alice colors $w$ red, Bob colors $u$ red, and we are back to a 2-coloring of $G$ and we can start the process over with Alice moving next.\\

Case 2) If Bob colors $u$ green, then Alice colors $w$ green and Bob colors $v$ blue. Again, we  are back to a 2-coloring of $G$ and we can start the process over with Alice moving next. $\Box$ \\

%

Of course, for all $G$, $\chi_g(G) \leq \chi_g^\infty(G)$. 

\begin{conjecture} \label{game-conj}For every graph $G$,  $\chi_g(G) < \chi_g^\infty(G)$.
\end{conjecture}




One might suspect that  $\chi^\infty(G) \leq \chi_g^\infty(G)$ for all $G$. However this is not always the case. From \cite{WFK} we know that $\chi^\infty(K_{1, 7}) = 4$, but below we show that $\chi_g^\infty(K_{1, 7}) > 4$. On the other hand, $\chi^\infty(P_3) = 4$  and $\chi_g^\infty(P_3) =3$.

Obviously $\chi^\infty(K_n) = \chi_g^\infty(K_n)=n+1$.  Similarly, $\chi^\infty(K_n-e) = \chi_g^\infty(K_n-e)=n+1$. 

We sometimes call the previous game the A-game, since Alice starts. We can also define the following variations on the A-game:\\

\noindent  $\bullet$ the B-game: same as the A-game except that Bob goes first.\\
\noindent $\bullet$ the A'-game: the same as the A-game except that Alice goes first in each round (so Alice may have two consecutive turns in the case $n$ is odd).\\
\noindent $\bullet$  the B'-game: same as the A'-game except that Bob goes first in each round.\\

Note that $P_3$ is an example in which the A-game and B'-game differ, since $\chi_g^\infty(P_3)=3$, but four colors are needed if Bob goes first. To see the latter, let $P_3=uvw$ and suppose Alice can win with three colors. Bob colors $v$ red, Alice colors $u$ green, and Bob colors $w$ blue. Bob then attempts to re-color $v$ and but cannot properly color it with the available colors, resulting in a win for Alice.

In fact,  $P_3$ is an example in which the A-game and B-game differ.  Let $P_3=uvw$ and suppose Alice can win the B-game with three colors. Bob colors $v$ red, Alice colors $u$ green, and Bob colors $w$ blue. Alice must now color $w$ green (same as coloring $u$ blue) else she loses the game. Bob now colors $u$ blue, and Alice loses the game on her next turn.

 However in general, we suspect there is little difference -- perhaps one additional color in some cases -- in the various games defined in this section. In the next section, we consider the A-game. In Section \ref{alt} we consider variations on the game and in Section \ref{quest}, we propose two more.

%
%
%
%

\section{Results}

\begin{theorem}\label{thm1}
Let $G$ be a graph with maximum degree $\Delta(G)$. Then $\chi_g^\infty(G) \leq \Delta(G)+2$.
\end{theorem}

\noindent{\it Proof:}
The standard proof technique showing that $\chi(G) \leq \Delta(G)+1$ can be applied to show that  $\chi_g^\infty(G) \leq \Delta(G)+2$.
$\Box$

\begin{proposition}\label{path}
Let $P_n$ be a path with $n \geq 4$. Then $\chi_g^\infty(P_n) =4$.
\end{proposition}

\noindent{\it Proof:}
That four colors suffice follows from Theorem \ref{thm1}. To see that four colors are  necessary, suppose to the contrary we can color $P_4$ with three colors; the argument is similar when $n > 4$ (and, in fact, easier in the odd cases). Let $P_4 = uvwx$. There are two cases depending on which vertex Alice colors first in round 1 (up to symmetry). If she colors $u$ with, say, color 1, then Bob colors $w$ with color 2, Alice colors $x$ with color 1 (or color 3, but color 1 is a better choice), and Bob colors $v$ with color 2. On the other hand, if Alice colors $v$ with color 1 first, then Bob colors $x$ with color 2, then Alice colors $u$ with color 3 (or 2, it doesn't matter), and Bob color $w$ with color 3. In either case, there is a vertex with all three colors in its closed neighborhood. Let us suppose that vertex is $w$, Alice must then change the color of a neighbor of $w$. But then Bob changes the color of the other vertex in $w$'s neighborhood. This  then will force Bob to win when $w$ is re-colored. $\Box$

\medskip

Any cycle is an example with $\chi_g^\infty(G) = \Delta(G)+2$. One can verify that $\chi_g^\infty(P_n)=4$ when $n > 3$. For comparison sake, it is known that for every graph $G$, $\chi (G)\leq \chi _{g}(G)\leq \Delta (G)+1$, see \cite{faig}.
We ask whether or not it is always true that for any cubic graph $G$,  $\chi_g^\infty(G)=5$?

The graph $K_{1, n}$ is called a {\it star}. When $n > 1$, the unique vertex of degree greater than one is called the \emph{center}.
We note that $\chi^\infty(K_{1, n})=4$ when $n>2$ (see \cite{WFK}) and $\chi_g(K_{1,n})=2$  when $n \geq 1$ (see \cite{des}).  The eternal game chromatic number of stars, however, is quite different, as we show next. Recall that above we showed that 
$\chi_g^\infty(K_{1, 2})=3$.

\begin{theorem}\label{lb}
When $n$ is odd $\chi_g^\infty(K_{1, n}) \geq \lceil \frac{n}{2} \rceil  +2$.
When $n \geq 4$ is even $\chi_g^\infty(K_{1, n}) \geq \frac{n}{2} +3$.
\end{theorem}

\noindent{\it Proof:}
We first explain the case when $n$ is even. Bob's strategy in round 1 is to color at least $n/2$ of the leaves using distinct colors for each -- or as many distinct colors as possible depending on the number of available colors (note that Alice may try to use one or more of these colors on other leaves). He then will force an additional color by choosing the center vertex on the first turn in round 2.  It is easy to see that Bob will win the game unless at least  $\frac{n}{2} +2$ colors are available. We can do one better than that by a more careful analysis which we now perform which will show that $\frac{n}{2} +3$ colors are needed. Suppose to the contrary that only
$\frac{n}{2} +2$ colors are available.

There are two cases.

\noindent{\bf{Case 1.}}  Suppose during round 1, Alice does not color the center vertex. Then Bob and Alice use $\lceil \frac{n}{2} \rceil+ 1$ colors on the leaves, a different color on the center vertex, and then on Bob's first move during round 2, he picks the center vertex and forces
an additional color.\\

\noindent{\bf{Case 2.}}  Suppose during round 1, Alice does color the center vertex (and we can assume without loss of generality that she does so on her first move). If Alice colors all leaves the same color, at the end of round 1, at most  $\lceil \frac{n}{2} \rceil+ 1$ colors are used.
 Bob starts round 2 and plays a new color on the center vertex. He can then ensure that all $\lceil \frac{n}{2} \rceil+ 2$ are used on vertices by the end of round 2. In round 3,  Alice cannot play the center vertex first, else an extra color is needed, so she plays a leaf, changing it from $c$ to $c'$. Bob then executes a ``mirroring" strategy, changing a leaf with color $c'$ to $c$. At the end of the round,  all $\lceil \frac{n}{2} \rceil+ 2$ are used on vertices. Bob then forces a new color on the first move of round 4.

\medskip

Now suppose $n$ is odd and thus Alice will go first on each round. We are claiming in this case that $\chi_g^\infty(K_{1, n}) \geq \lceil \frac{n}{2} \rceil  +2$. Suppose to the contrary that only $\lceil \frac{n}{2} \rceil+ 1$ colors are available. 

Bob will again try to use as many different colors as possible on the leaves.
There are two cases.

\noindent{\bf{Case 1.}}  Suppose Bob colors the leaves with $\frac{n}{2} $ distinct colors in round 1 (which occurs if Alice colors the central vertex in round 1). Then the total number of colors used in round 1 is  $\frac{n}{2}  +1$. If Alice colors the center vertex on the first move
in round 2, then an additional color is needed and we are done. So suppose Alice colors a leaf on her first move in round 2. Whatever move she makes, Bob mirrors it on his next move (e.g., if she re-colors a vertex from $c$ to $c'$, Bob will re-color a vertex from $c'$ to $c$.). Eventually, Bob can re-color the center vertex on his last turn in round 2, forcing the use of an additional color.

\noindent{\bf{Case 2.}}  Suppose Bob colors the leaves with $\frac{n}{2}  -1$  distinct colors in round 1 (which occurs if Alice does not color the central vertex in round 1 and she uses at least one of the same colors as Bob on one or more leaves). Then each of these colors Bob
uses is distinct from the color Alice uses on the first leaf. Bob then uses color $\frac{n}{2}  +1$ on the last move of round 1. The remainder proceeds as in Case 1.
$\Box$

\medskip

We next show the lower bound from the previous result is tight.

\begin{theorem}\label{star-lower}
When $n$ is odd $\chi_g^\infty(K_{1, n}) \leq \lceil \frac{n}{2} \rceil  +2$.
When $n \geq 4$ is even $\chi_g^\infty(K_{1, n}) \leq \frac{n}{2} +3$.
\end{theorem}

\noindent{\it Proof:} We assume $n > 2$ as the result is easy to see when $n \leq 2$. 
Let $c(v)$ denote the color assigned to $v$ at any given time (initially null) and in each case we shall assume the number of
colors specified in the statement of the theorem are available.

First suppose $n$ is odd, which means Alice goes first on each round. Suppose without loss of generality that the
coloring of the graph done in round 1 uses $\lceil \frac{n}{2} \rceil+1$ colors, because Alice's strategy
is to color leaves with color 1 on each of her turns. Thus, $c(v) = 1$
for $\lceil \frac{n}{2} \rceil$ of the leaves. We may assume that Bob uses a distinct color other than 1 on each of his turns on each of the
remaining vertices. Alice's strategy in the next round is to change any vertex not colored 1 
to 1, and then Bob will do one of two things on his first turn of round 2.

\noindent{\bf{Case 1)}} If Bob decides to choose the center vertex, he must choose a $\lceil \frac{n}{2} \rceil+2^{nd}$
color. Then, Alice chooses another vertex $v$ with $c(v) \neq 1$ and re-colors it to 1, and Bob (in order to maximize the number of colors used) either (a) mirrors
her actions by choosing a vertex $v$ with $c(v) = 1$ and changing it to the color (which is not equal to 1)  that was
the previous color on the vertex Alice chose or (b) he could also change it to the
previous color of the center vertex. He continues to ensure that he does not use
the same color twice during the round. When Alice runs out of vertices that are not colored 1, she chooses one of the remaining two vertices (both of color 1) to change to any of the colors Bob used on another leaf during
round 2. On Bob's final turn of the round, there are two colors of the $\lceil \frac{n}{2} \rceil + 2$ used thus far
that have not yet been used during round 2. This is because Bob had one fewer
turn on the leaves than Alice (so he could not counteract her actions
each time), and there is already one more color in the set than necessary for
the initial coloring. Thus, Bob picks one of the remaining colors, and the other
is free for the same process to repeat in the following round without needing to
force the use of any more colors. Therefore, $\lceil \frac{n}{2} \rceil + 2$ colors suffice.\\

\noindent  {\bf{Case 2)}} If Bob selects a leaf to re-color, then in order to maximize the number of colors used on the leaves,
 he will use a mirroring strategy to counteract Alice's turns (i.e, if Alice re-colors a vertex from $c$ to $c'$ , he will do the opposite). This can continue on the leaves 
until one remains of color 1 and it is Alice's turn. She will pick a color
which Bob used to re-color another leaf. Bob then chooses the center vertex and
uses a $\lceil \frac{n}{2} \rceil+2^{nd}$ color. This leaves the initial color of the center vertex free
for the same process to repeat on the next round. $\lceil \frac{n}{2} \rceil + 2$ colors suffice.\\

\medskip

Now suppose $n$ is even, which means Alice goes first on each odd-numbered round.
Again, Alice's general strategy is to color as many leaves as possible with the same color.
Since she can color at least $\frac{n}{2}$ of them the same color, at most $\frac{n}{2} +1$ different colors are used
on the leaves at the end of any given round, plus possibly an additional color for the center vertex.  Since the center vertex needs to change colors in each round,
we claim that $\frac{n}{2} +3$  colors suffice. To see this, note that at some point in a round in which Bob goes first, Bob can change the color of a leaf, thus (temporarily) there may be
$\frac{n}{2} +2$ different colors on leaves. However, on Alice's next move, she can choose a leaf to re-color that reduces the number of different colors on the leaves to $\frac{n}{2} +1$.
$\Box$

\medskip

It is easy to verify that $\chi_g^\infty(K_{2, 2})=4$.\\


It is known that $\chi_g(K_{n, n})=3$ when $n$ and $m$ are both larger than 1, see \cite{des}.

\begin{proposition}\label{comp}
$\chi_g^\infty(K_{n, n}) \leq 5$ when $n \geq 5$. 
\end{proposition}

\noindent{\it Proof:} Let the vertex parts of $K_{n, n}$ be $A=\{a_1, \ldots, a_n\}$ and $B=\{b_1, \ldots, b_n\}$. 
Alice's strategy is to ensure that at the end of each round, there are at least two distinct colors used on $A$ and at least two distinct colors used on $B$. This then allows vertices in each part to switch back and forth between those two colors each time they are re-colored.
This is easily seen to be possible when $n \geq 5$. We note that in total, five colors are needed, since Bob can force a third color onto either $A$ or $B$ whilst Alice is trying to ensure there are at least two colors on both $A$ and $B$.
$\Box$

\medskip

This demonstrates that $\chi_g^\infty$ does not have hereditary properties, since $K_{1, n}$ is an induced subgraph of $K_{n, n}$.

\begin{proposition}
$\chi_g^\infty(K_{n, n}) \geq 4$ when $n \geq 3$. 
\end{proposition}

\noindent{\it Proof}: Suppose to the contrary that three colors suffice. Let the vertex parts be $A=\{a_1, \ldots, a_n\}$ and $B=\{b_1, \ldots, b_n\}$. 
Observe that Bob can ensure that at the end of round 1, vertices in $A$ use two distinct colors -- colors 1 and 2 -- and each vertex in $B$ is color 3 (without loss of generality). On Alice's first move in round 2, she cannot re-color any vertex in $B$, as this would require a fourth color.
If she re-colors a vertex in $A$ so that there are still two distinct colors in $A$, Bob will be unable to move on his turn if he chooses a vertex in $B$. So suppose Alice re-colors, say, $a_1$ so that all the vertices in $A$ are the same color, say color 2. Bob now re-colors a vertex in $B$ with color 1, so that $B$ now contains
vertices of two distinct colors. Alice now needs a fourth color if she chooses a vertex in $A$, so she chooses a vertex in $B$. Bob now chooses a vertex in $A$ and needs a fourth color, winning the game.

We note that if $n \geq  4$,  Bob can ensure that at the end of round 1, vertices in $A$ use two distinct colors 1 and 2 and vertices in $B$ use two distinct colors. Therefore, Bob can win the game more quickly in round 2. We discuss this issue further in the last section of the paper.
$\Box$

\begin{proposition}
$\chi_g^\infty(K_{n, n}) \geq 5$ when $n \geq 5$. 
\end{proposition}

\noindent{\it Proof}: Suppose to the contrary that four colors suffice. Let the vertex parts be $A=\{a_1, \ldots, a_n\}$ and $B=\{b_1, \ldots, b_n\}$. Let the color set be $1, 2, 3, 4$ and let $c(v)$ denote the color assigned to $v$ at any given time (initially null). First suppose that at the end of round 1,
there is only one color used on $A$, say color 1. In this case, Bob can ensure that at least three colors are used on $B$. If Alice chooses a vertex in $A$ for her first move in round 2, a fifth color must be used. Hence she must choose a vertex in $B$. In other to avoid a fifth color being played on $A$ on Bob's first move in round 2, Alice must change the color of a vertex $v \in B$ such that $v$ is the only vertex in $B$ with that color, say color 4, which she changes to 3. Note that Bob can have assured that there are at least two vertices in $B$ of color 2 and at least two vertices of color 3 at the end of round 1. 
Bob now changes a vertex in $A$ to color 4, and Alice must use a fifth color on her next turn.

We now claim that in fact Bob can force there to be only one color on $B$. Assume without loss of generality that Alice colors $a_1$ first with color 1. Then Bob chooses $a_2$ and colors it 2. If  Alice chooses $b_1$ and colors it 3, and then Bob colors $a_3$ with color 4, then all vertices in $B$ will have to be colored with 3 in round 1. If, on the other hand, Alice choose $a_3$ and colors it 1 or 2 (note that if she colored it with color 3, then Bob would color $a_4$ with 4 and there would be no colors available for $B$), then Bob colors $a_4$ with 4, again forcing all of $B$ to be the same color. 
$\Box$

\begin{proposition}\label{k44}
$\chi_g^\infty(K_{4, 4}) = 4$. 
\end{proposition}

\noindent{\it Proof}:  We show that $\chi_g^\infty(K_{4, 4}) \leq 4$. If Alice can force there to be at least two distinct colors used on $A$ and two distinct colors used on $B$ at the end of round 1, she can successfully proceed with 4 colors, as described above. Suppose this is not the case and that only 1 color, say color 1, is used on $A$
during round 1. In order for Bob to do this, Alice can ensure that only colors 2 and 3 are used on $B$ during round 1. On her first move in round 2, Alice re-colors a vertex in $A$ with color 4. This then forces the vertices in $B$ to swap between colors 2 and 4 and color 4 to be used on the remainder of $A$. The process can be repeated in subsequent rounds.
$\Box$

\begin{proposition}
$\chi_g^\infty(K_{3, 3}) = 4$. 
\end{proposition}

\noindent{\it Proof}:  Similar to Proposition \ref{k44}.
$\Box$

\medskip


\begin{theorem}
Let $G$ be a connected graph other than $K_1, K_2$, or $P_3$. Then $\chi_g^\infty(G) >  3$.  
\end{theorem}

\noindent{\it Proof:} We can assume $G$ has at least four vertices, since $\chi_g^\infty(K_3) > 3$. We also assume that $\Delta(G) > 2$, else $G$ is a cycle and thus four colors are needed. Of course, $G$ may be a tree, but not necessarily so. Let $v$ be a vertex of degree at least three, with $x, y, z$ three of the neighbors of $v$. Suppose to the contrary that 3 colors suffice. Our general strategy is as follows. Assume without loss of generality that $c(v)=1$. Bob can force at least two distinct colors, other than color 1, be used on $\{x, y, z\}$ in round 1. If Bob moves first in round 2, he can force a fourth color by choosing $v$. Otherwise, Bob can still force that, at any time, it is either Alice's turn and there are two distinct colors on $\{x, y, z\}$ or it is his turn and there are either (i) two distinct colors on $\{x, y, z\}$  (ii) there is only one distinct color on $\{x, y, z\}$ but at least one of these three vertices has yet to be re-colored. Hence a fourth color will be needed unless there is a vertex $w \neq v$ that is also adjacent to each of $x, y, z$ and Alice manages to color $w$ with color 2, whilst $c(v)=1$. This then forces $x, y, z$ to all have color 3. If Bob gets to color one of these five vertices first in round 2, obviously he forces a fourth color. So assume Alice colors one of these first. Then it must be that she re-colors $v$ with color 2 (or equivalently $w$ with color 1). But then Bob can re-color $w$ with color 1 (equivalently $w$ with color 2) and then a fourth color will be forced on one of $x, y, z$.
$\Box$

\begin{proposition}
(a) $\chi_g^\infty(K_{2, n}) = 4$ if $n \in \{2, 3\}$.\\
(b) $\chi_g^\infty(K_{2, n}) \geq n/2 + 3$ if $n > 3$ is even.\\
(c)   $\chi_g^\infty(K_{2, n}) \geq \lfloor \frac{n}{2} \rfloor + 3$ if $n > 3$ is odd.
\end{proposition}

\noindent{\it Proof}:  Part (a) is not difficult to verify by hand. For parts (b) and (c), Let the vertex parts be $A=\{a_1, a_2\}$ and $B=\{b_1, \ldots, b_n\}$. Let us first suppose $n$ is even. Bob's basic strategy in both (b) and (c) will be to make sure both vertices in $A$ have the same color and then use as many different colors on $B$ as possible (which will be $n/2 + 1$).

For part (b), observe that Bob can  ensure both vertices in $A$ have the same color after round 1 (and, in fact, after all subsequent rounds). Therefore, Bob can ensure $n/2 + 2$ colors are used in round 1 in total ($n/2 + 1$ on $B$ and one on $A$). An additional color will be needed on $A$ in round 2, regardless of when a vertex in $A$ in chosen (since if Alice chooses to play a vertex $v \in B$ first in round 2, she must change its color to that of another vertex $u \in B$; then on Bob's turn he can change $u$'s color to the original color of $v$).

For part (c), the argument is similar to (b), except that Bob cannot necessarily force both vertices in $A$ to be the same color after each round (otherwise, Bob can force $\lfloor \frac{n}{2} \rfloor + 2$ colors be used in round 1 and an additional color in round 2). So suppose the vertices in $A$ have two different colors at the end of some round. One can easily verify in this case that Bob can ensure that at least $\lfloor \frac{n}{2} \rfloor + 1$ colors are used on $B$.
$\Box$

\section{Alternate Games}\label{alt}

An alternate version of game chromatic number called {\it game chromatic number II} was introduced in \cite{chen}. In this version,  Bob can only use the colors that have been already used on the graph, unless he is forced to use a new color to
 guarantee that the graph is colored properly. The number of colors needed for $G$ in this game is denoted by  $\chi_g^{*}(G)$.  It was shown in \cite{chen} that $\chi_g^{*}(T) \leq 3$ for any tree $T$, as compared to the bound of 4 for $\chi_g(T)$, see \cite{faig}.
 
We define the {\it eternal vertex coloring game II} in the obvious manner and use $\chi_g^{\infty*}(G)$ to denote the number of colors needed when this game is played eternally on $G$: that is, the game is played over a series of rounds and in each round Bob can only use colors that have already been introduced unless he is forced to use a new color. Of course,  $\chi_g^{\infty*}(G) \leq \chi_g^{\infty}(G)$.
$K_3$ is an example where  $\chi_g^{\infty*}(G) < 2\chi_g^*(G)$ and $P_4$ is an example where  $\chi_g^{\infty*}(P_4) = 4 = 2\chi_g^*(P_4)$.  

We shall consider stars again and begin by noting that $\chi_g^{\infty*}(K_{1, 3})=3$ and $\chi_g^{\infty*}(K_{1, 4})=4$.  
Interestingly, we shall see in next result that it may take Bob several rounds to force the number of colors needed. The next result is also a good example of the power Bob gets when the number of vertices is odd.

\begin{theorem}
Then (a) if $n$ is odd then $\chi_g^{\infty*}(K_{1, n})=3$ and (b) if $n \geq 4$ is even then $\chi_g^{\infty*}(K_{1, n}) =  \frac{n}{2}+2$.
\end{theorem}

[Note: does this proof need work?]

\noindent{\it Proof:} 
Part (a) is easy to see, noting that Alice can choose the center vertex of the star on the first turn in each round, coloring it with the smallest color possible. For part(b), note that two colors are used in round 1. Bob can then force colors 3 and 4 to be used during round 2: he does so by choosing a leaf on his first turn and the center vertex on his second turn. Subsequent to that, for the next several rounds, Bob can force an additional color to be used on either each round or every other round by coloring the leaves with as many different colors as possible during each round, eventually forcing the center vertex to use a new color -- it may need to wait until a round in which Bob goes first for each new color to be introduced (as in round 2).
Bob can eventually force the use of $\frac{n}{2}$ different colors on the leaves, once that many colors have brought into play. At which point the game is the same as the A-game and thus the result follows as in Theorem \ref{lb}.
$\Box$

\medskip

One could also consider a more restrictive version of this game in which Bob must color whatever vertex he chooses with the smallest color possible (as opposed to being able to use any color, once a color becomes allowable).  We call this the {\it Greedy Coloring Game}.
Denote the number of colors needed for this game as $\chi^2_g(G) $ when only one round is played and $\chi_g^{\infty2}(G)$ when the game is played eternally.

It is then natural to ask if there exist any graphs 
where this restriction changes the number of colors needed versus $\chi_g^{\infty*}(G)$. The answer is ``yes", as we show next.

\begin{theorem}
$\chi_g^{\infty2}(K_{1, n})=3$ when $n$ is odd and $\chi_g^{\infty2}(K_{1, n})=4$ when $n \geq 4$ is even.
\end{theorem}

\noindent{\it Proof:} 
First consider $K_{1, n}$ when $n > 3$ is odd and thus Alice goes first in each round. In this case, only 3 colors are needed. Alice chooses the center vertex of the star first in each round, giving it color 1 in odd numbered rounds and all other vertices getting color 2) and color 3 in even numbered rounds (and all other vertices getting color 1).

When $n \geq 4$ is even, we show $K_{1, n}=4$. That fours colors are necessary is easy  to see, so we show that four colors suffice. On round 1, Alice plays color 1 on the center and then color 2 gets played on all the leaves. In round 2, Bob goes first. If he plays the center, he will use color 3 and then color 1 gets played on all the leaves (which is no advantage to Bob). So assume Bob plays a 3 on a leaf. Then Alice plays 1 on a leaf, Bob plays 4 on the center, and 1's gets played on the remaining leaves. Alice, starting the next round, re-colors the lead that is color 3 with color 2. Then either remaining leaves get color 2 and the center gets color 1 (which puts us back in a previous configuration), or Bob plays 3 on the center on some turn in this round and the remainder of the leaves get color 1. However, this coloring completes the round with only two different colors, 1 and 3, on the vertices, which is essentially the same at this point as the configuration in which only 1 and 2 are used.
$\Box$

\begin{conjecture} For every graph $G$,  $\chi^2_g(G) < \chi_g^{\infty2}(G)$.
\end{conjecture}
%



One could further restrict the rules so as to force Alice to also choose smallest color for each vertex she chooses.  We call this the {\it Very Greedy Coloring Game}. Denote the number of colors needed for this game as $\chi^3_g(G)$ when only one round is played and $\chi_g^{\infty3 }(G)$ when the game is played eternally. 
Recall that the {\it coloring number} of $G$, denoted $col(G)$, is the smallest integer $k$ such that every subgraph of $G$ has a vertex of degree less than $k$. 
Though intuitively, there seems to be some relationship between $col(G)$ and $\chi^3_g(G)$, in fact they are not the same:  $col(K_{n, n})=n+1$ and $\chi^3_g(K_{n, n})=2$; whereas  $col(P_4)=2$ and $\chi^3_g(P_4)=3$.

It is clear that $\chi_g^{\infty2}(G) \leq \chi_g^{\infty3}(G)$ and $\chi_g^{\infty2}(G) \leq  \chi_g^{\infty*}(G) \leq \chi_g^{\infty}(G)$ for all $G$. 

\begin{question}
Is it true that $\chi_g^{\infty3}(G) \leq  \chi_g^{\infty*}(G)$ for all graphs $G$?
\end{question}

\begin{theorem}
When $n$ is odd, $\chi_g^{\infty3 }(K_{1, n})=3$, when $n \geq 4$ is even, $\chi_g^{\infty3 }(K_{1, n})=4$.
\end{theorem}

\noindent{\it Proof:} The case when $n$ is odd is trivial. So suppose $n \geq 4$ is even. It is easy to see that $\chi_g^{\infty3 }(K_{1, n}) \geq 4$, as Bob can force color 4 to be used during round 2. In order to see that four colors suffice, the key observation is that Alice can maintain the invariant that, at the end of each round and prior to the re-coloring of the center vertex during each round, the leaves are colored with at most two distinct colors.  We strengthen this invariant as follows:  (a) no leaf is ever color 4 and (b) if the center is color 1 after a round and the leaves have two distinct colors, then the color that appears on fewer leaves appears at most twice.

Note that in order for Bob to force a fifth color, it would have to be the case that there are three distinct colors on the leaves. Of course, after the center vertex is re-colored during a round, it may be possible that three distinct colors are used on the leaves for a short time, until the end of the round at which point we claim that our invariant will be true. We now prove that the invariant can be maintained, from which the theorem follows. It is trivially true during and after round 1: Alice ensures  this by coloring the center first and then all the leaves must get color 2. Now let us move past round 1. If all the leaves are the same color after some round, then they all get changed to the same color during the next round until the center vertex is re-colored, at which point the remaining leaves may get changed to a second color -- however, by Alice choosing the center vertex as soon as possible, our invariant is assured. We note that if there are two distinct colors on the leaves at the end of the round in this case, one of these colors on the leaves must be color 1.

Now suppose the leaves are colored with two colors after some round. It is easy to see that we never use color 4 on a leaf: since when assigning the color to a leaf, one only needs to consider its current color and the color of the center; hence one of colors 1, 2, 3, will be available to assign to each leaf.

First suppose the center vertex is some color $c > 1$ and some leaf is color 1. If the center is re-colored first (which is Alice's preference if she goes first), it cannot get color 1 (since some of the leaves have color 1), which means that each leaf will be re-colored color 1. 
If the center vertex is not re-colored first (meaning Bob goes first), then the leaves that are re-colored prior to the center getting re-colored get color 1, the center is re-colored some color other than 1, and the remaining leaves get color 1. Thus the invariant is maintained.

Now suppose the center vertex is some color $c > 1$ and no leaf is color 1. Then it must be the case that either all the leaves have the same color, which we addressed earlier, or else some leaf is color 4, which we know cannot happen.


Next suppose the center vertex is color 1.  If the center is re-colored first (which is Alice's preference if she goes first), it gets some color which is either 2, 3, 4 (since the leaves use at most two different colors). Therefore the leaves all get color 1 during this round. If the center is not re-colored first (meaning Bob goes first), there are two cases.\\


\noindent {\bf{Case 1.}} Suppose Bob first chooses a color that appears on at most two leaves. Then Alice simply ``mirrors'' Bob's moves during the round, choosing another leaf of the same color Bob just chose. This guarantees the invariant is maintained.\\

\noindent {\bf{Case 2.}} Suppose Bob first chooses a color that appears on more than two leaves. Note that this may create a situation where there are three distinct colors on the leaves -- otherwise Bob has changed the color of a leaf to the color that exists on another leaf and Alice can simply chooses the center vertex on her next move and we will be guaranteed the invariant is maintained. (Or alternatively, Alice can essentially duplicate what Bob just did, choosing a leaf of the same color Bob just chose and this will also allow the invariant to be preserved). So suppose we now have three colors on the leaves -- the only way this occurs is if the leaves were colored with colors 2 and 4 before the round started. Since no leaves are color, 4, we are done.
$\Box$

\begin{question}
Can we differentiate between  $\chi_g^{\infty2}(G)$ and $\chi_g^{\infty3 }(G)$? That is,  does there exist a graph $G$ such that $\chi_g^{\infty2}(G) \neq \chi_g^{\infty3}(G)$?
\end{question}

We observe that $\chi^3_g(T) \leq 3$ for any tree $T$, as the proof that $\chi_g^{*}(T) \leq 3$ applies to $\chi^3_g(T)$ (and thus also for $\chi^2_g(T)$).

\begin{question}
Is there an integer $c$ such that for all trees $T$ with $\chi_g^{\infty3}(T)  \leq c$? Likewise $\chi_g^{\infty2}(T)$.
\end{question}


%



Caterpillars are one class of trees that are easy to analyze, as it is not hard to show that $\chi_g^{\infty3}(T)  \leq 6$ for any caterpillar $T$.

\begin{conjecture} (a) For every graph $G$,  $\chi^3_g(G) < \chi_g^{\infty3}(G)$. (b)  For every graph $G$,  $\chi^2_g(G) < \chi_g^{\infty2}(G)$. (c)  For every graph $G$,  $\chi^*_g(G) < \chi_g^{\infty*}(G)$.
\end{conjecture}

\begin{problem}
 Characterize for each of these models discussed in this section the graphs needing three or four colors.
\end{problem}

\begin{question}
Let $G$ be a graph with subgraph, or induced subgraph, $H$. Is it necessarily true that $\chi_g^{\infty2}(H) \leq \chi_g^{\infty2}(G)$? Is it necessarily true that $\chi_g^{\infty3}(H) \leq \chi_g^{\infty3}(G)$ ?
\end{question}

\section{Further Questions}\label{quest}

\begin{question}
What is  the computational complexity of each of the eternal coloring games?
\end{question}

\begin{question}
Characterize the connected graphs (or the trees) with  $\chi_g^\infty(G) =  4$. 
\end{question}


In seems an interesting future direction to consider $\chi_g^{\infty2}$, $\chi_g^{\infty3}$, and $\chi_g^{\infty}$ as well as $\chi_g^{2}$, $\chi_g^{3}$, and $\chi_g^{*}$ for various classes of graphs. For example,
it is easy to show that for a grid graph $G$ with a sufficiently large number of rows and  columns,  $5 \leq \chi_g^\infty(G) \leq 6$.  

\begin{problem}
Determine a tight bound on  $\chi_g^{\infty2}$, $\chi_g^{\infty3}$, and $\chi_g^{\infty}$ for large grid graphs. We conjecture at most five colors suffice for the greedy game and probably also for the very greedy game.
\end{problem}

\begin{question} Is it true that if $\chi_g^\infty(G)=k$ then Alice can win the A-game on $G$ with a color set of $k+1$ colors?
\end{question}

In the following question, we note the proof of Theorem \ref{lb}, in which it appears to take Bob four rounds to win on a star.

\begin{question} Suppose the A-game is played on some graph $G$ with $k < \chi_g^\infty(G)$ colors and Bob wins the game. Is it necessarily the case that Bob can win the game prior to the beginning of round 3? 
What is the maximum number of rounds it  takes Bob to win to any graph?
\end{question}

\begin{question}
Are there any graphs $G$ with $\chi^\infty(G) = \chi_g^\infty(G)$ and $\chi^\infty(G)  < \Delta(G)+2$?
\end{question}

\begin{question}
Can we characterize the graphs with $\chi(G) + 1 = \chi_g^\infty(G)$? Likewise for the other graph coloring games introduced in this paper.
\end{question}


We raise an issue regarding the definition of the eternal vertex coloring game (and some of its variants). Suppose Bob chooses vertex $v$ to re-color (in round 2 or later) and Bob has no legal move on this vertex, however Bob does have a legal move on vertex $u \neq v$.  As we have defined the game now, Bob wins. But should we force Bob to choose some other vertex, such as $u$, where he has a legal move? It may be worth exploring this version of the game. One way to do this is to consider the alternate game discussed in Section \ref{alt}. Two obvious other variations come to mind:

$\bullet$ (Strong Eternal Game Coloring) We require that Bob must choose a vertex on his turn that can be re-colored properly, if one exists (as opposed to being able to choose any vertex); or

$\bullet$ (Ordered Eternal Game Coloring) The vertices must be chosen in the same  order (pre-determined or fixed after round 1) each round. 

We know that  $\chi^\infty(C_5)=\chi_g^\infty(C_5)=4$. However, if one follows the Strong Eternal Game Coloring rules, three colors suffice. To see this, a 3-coloring results from round 1. On the first move in round 2, Bob must choose a vertex that can be re-colored (a vertex adjacent to two vertices of the same color), as opposed to a vertex adjacent to two different colors (which would have necessitated the use of a fourth color). It is not hard to see that round 2 can be successfully completed with three colors and the subsequent rounds can continue in a similar manner.


%

\end{document}